\documentclass[11pt,reqno,tbtags,a4paper]{amsart}
\usepackage{amssymb}
\usepackage{xpunctuate}
\usepackage{url}
\usepackage[square,numbers]{natbib}
\bibpunct[, ]{[}{]}{;}{n}{,}{,}

\title
{On the independence number of some random trees}

\date{19 March, 2020; revised 20 March, 2020}

\author{Svante Janson}
\thanks{Partly supported by the Knut and Alice Wallenberg Foundation}
\address{Department of Mathematics, Uppsala University, PO Box 480,
SE-751~06 Uppsala, Sweden}
\email{svante.janson@math.uu.se}
\newcommand\urladdrx[1]{{\urladdr{\def~{{\tiny$\sim$}}#1}}}
\urladdrx{http://www2.math.uu.se/~svante/}


\numberwithin{equation}{section}

\renewcommand\le{\leqslant}
\renewcommand\ge{\geqslant}

\allowdisplaybreaks

\theoremstyle{plain}
\newtheorem{theorem}{Theorem}[section]

\theoremstyle{definition}

\newtheorem{exampleqqq}[theorem]{Example}
\newenvironment{example}{\begin{exampleqqq}}
  {\hfill\qedsymbol\end{exampleqqq}}

\newtheorem{remarkqqq}[theorem]{Remark}
\newenvironment{remark}{\begin{remarkqqq}}
  {\hfill\qedsymbol\end{remarkqqq}}

\newtheorem*{ack}{Acknowledgement}

\theoremstyle{remark}

\newenvironment{romenumerate}[1][-10pt]{
\addtolength{\leftmargini}{#1}\begin{enumerate}
 \renewcommand{\labelenumi}{\textup{(\roman{enumi})}}%
 }{\end{enumerate}}

\newenvironment{PXenumerate}[1]{
\begin{enumerate}
 \renewcommand{\labelenumi}{\textup{(#1\arabic{enumi})}}%
 }{\end{enumerate}}

\newcounter{oldenumi}
{\setcounter{oldenumi}{\value{enumi}}
\begin{romenumerate} \setcounter{enumi}{\value{oldenumi}}}
{\end{romenumerate}}

\newcounter{thmenumerate}

\newcounter{xenumerate}   

\newcommand{\refT}[1]{Theorem~\ref{#1}}

\newcommand{\refR}[1]{Remark~\ref{#1}}

\newcommand{\refS}[1]{Section~\ref{#1}}
\newcommand{\refSs}[1]{Sections~\ref{#1}}

\begingroup
  \divide\count255 by 60
  \multiply\count255 by -60
  \advance\count255 by \time
  \ifnum \count255 < 10 \xdef\klockan{\the\count1.0\the\count255}
  \else\xdef\klockan{\the\count1.\the\count255}\fi
\endgroup

\newcommand{\sumko}{\sum_{k=0}^\infty}

\newcommand{\sumid}{\sum_{i=1}^d}

\newcommand{\prodim}{\prod_{i=1}^m}

\newcommand\set[1]{\ensuremath{\{#1\}}}

\newcommand\xpar[1]{(#1)}
\newcommand\bigpar[1]{\bigl(#1\bigr)}
\newcommand\Bigpar[1]{\Bigl(#1\Bigr)}
\newcommand\biggpar[1]{\biggl(#1\biggr)}

\newcommand\bigsqpar[1]{\bigl[#1\bigr]}
\newcommand\Bigsqpar[1]{\Bigl[#1\Bigr]}

\newcommand\xcpar[1]{\{#1\}}

\newcommand\bigabs[1]{\bigl\lvert#1\bigr\rvert}
\newcommand\Bigabs[1]{\Bigl\lvert#1\Bigr\rvert}

\def\rompar(#1){\textup(#1\textup)}    

\newcommand\Bigparfrac[2]{\Bigpar{\frac{#1}{#2}}}
\newcommand\biggparfrac[2]{\biggpar{\frac{#1}{#2}}}

\def\xexp(#1){e^{#1}}

\newcommand\ntoo{\ensuremath{{n\to\infty}}}

\newcommand\ttoo{\ensuremath{{t\to\infty}}}

\newcommand\punkt{\xperiod}    
\newcommand\iid{i.i.d\punkt}    
\newcommand\ie{i.e\punkt}
\newcommand\eg{e.g\punkt}
\newcommand\viz{viz\punkt}
\newcommand\cf{cf\punkt}
\newcommand{\as}{a.s\punkt}

\newcommand{\tend}{\longrightarrow}

\newcommand\asto{\overset{\mathrm{a.s.}}{\tend}}

\newcommand\bbZ{\mathbb Z}

\newcounter{CC}
\newcounter{cc}

\newcommand\E{\operatorname{\mathbb E{}}}
\renewcommand\P{\operatorname{\mathbb P{}}}

\newcommand\Exp{\operatorname{Exp}}
\newcommand\Po{\operatorname{Po}}

\newcommand\ga{\alpha}
\newcommand\gb{\beta}
\newcommand\gd{\delta}
\newcommand\gD{\Delta}

\newcommand\gam{\gamma}

\newcommand\gl{\lambda}

\newcommand\gs{\sigma}

\newcommand\gth{\theta}

\newcommand\cT{{\mathcal T}}

\newcommand\indic[1]{\boldsymbol1\xcpar{#1}}

\newcommand\qw{^{-1}}
\newcommand\qww{^{-2}}
\newcommand\qq{^{1/2}}

\newcommand\intoi{\int_0^1}
\newcommand\intoo{\int_0^\infty}

\newcommand\oi{\ensuremath{[0,1]}}
\newcommand\ooi{(0,1]}
\newcommand\ooo{[0,\infty)}

\newcommand\setoi{\set{0,1}}

\newcommand\dd{\,\mathrm{d}}
\newcommand\ddx{\mathrm{d}}

\newcommand\ctt{\cT_t}
\newcommand\cttau{\cT_\tau}
\newcommand\intot{\int_0^t}
\newcommand\chix{\chi'}
\newcommand\LXi{L_\Xi}
\newcommand\AAA{\ref{BPfirst}--\ref{BPmub}}
\newcommand\fX{\mathfrak X}

\newcommand\gdp{\gD p}

\begin{document}

\begin{abstract} 
We show that for many models of random trees,
the independence number divided by the size converges almost surely to a
constant as the
size grows to infinity; the trees that we consider include
random recursive
trees, binary and $m$-ary search trees, preferential attachment trees, and
others. The limiting constant is computed, analytically or numerically, for
several examples. 
The method is based on Crump--Mode--Jagers branching processes.
\end{abstract}

\maketitle


\section{Introduction}\label{S:intro}

The \emph{independence number} 
\ie, the maximum size of an independent set of nodes,
is a  quantity 
that has been studied for various models of random trees 
(and other random graphs, not considered here). 
In the present paper we consider rooted trees that can be
constructed as family trees of a Crump--Mode--Jagers branching process
stopped at a suitable stopping time;
this includes, for example,  
random recursive trees, preferential attachment trees, 
fragmentation trees,
binary search trees and $m$-ary search trees;
see 
\refS{SSCMJ} and \cite{SJ306} for details,
and the examples in \refSs{SRRT}--\ref{Slast} below.

We denote the independence number of $T$ by $I(T)$.
Our main result, \refT{T1}, gives a strong law of large numbers for $I(T)$; more
precisely,  it shows convergence almost surely (a.s.)\ of $I(T_n)/|T_n|$,
the fraction of nodes that belong to
a maximum independent set, for a sequence $T_n$ of random trees. 
The limit $\nu$ is a constant depending on the random tree model; the
theorem expresses 
this limit in terms of the solution $p(t)$ of the functional
equation \eqref{tl1}. We show in \refSs{SRRT}--\ref{Slast} how this equation
can be
solved and $\nu$ found explicitly (at least numerically) 
in some important examples,
\viz{} random recursive trees, binary search trees,  preferential
attachment trees, extended binary search tress and $m$-ary search trees (in
particular $m=3$).

Note that the cases of random recursive trees and binary search trees have
been studied before.
For random recursive trees, the expectation was found already by
\citet{MM1975,MM1978}.
More recently, both  \citet{Dadedzi} and \citet{FuchsEtal} prove 
(independently, and with different methods)
the weak version (\ie, convergence in probability) of
\eqref{t1} below for random recursive trees and binary search trees, 
with explicit $\nu$,
and also a much stronger central limit theorem. 
Nevertheless, we think that the present approach is of interest, since it is
quite general; moreover, it gives convergence a.s. 
Furthermore, although we only prove a law of large numbers in
the present paper, we hope that future development of our methods will also
lead to a central limit theorem.

Similar results have also been proved for other types of random trees.
 For simply generated trees, see \eg{}
\citet{MM1973, 
MM1977I} 
and 
\citet{BanderierKP}; 
for uniform unlabelled trees (rooted or unrooted), see
\citet{MM1977II}. 

\begin{remark}
  As is well known, for trees, several other quantities are determined by the
  independence number by linear relations, and our results thus immediately
  transfer to these quantities.
These include, for example:
\begin{romenumerate}
  
\item 
 The \emph{matching number}, \ie, the maximum size of a partial
matching. This equals $|T|-I(T)$,
\item 
The \emph{minimum size of a vertex cover}, \ie, of a vertex set that
contains at least one end-point of every edge. This equals the matching
number, \ie, $|T|-I(T)$.
\item 
The \emph{nullity}, \ie, the dimension of the kernel of the adjacency
matrix, or
the multiplicity of the eigenvalue 0 of the adjacency matrix.
This equals $2I(T)-|T|$.
(The results in \cite{Dadedzi} referred to below are actually stated for the nullity.)
\end{romenumerate}
See \eg{} \cite{FuchsEtal} for further examples.
\end{remark}

\section{Preliminaries}\label{Sprel}

We give some definitions and notation, together with some 
known results  that will be used.

The number of nodes of a tree $T$ is denoted $|T|$.


If $T$ is a rooted tree, and $v\in T$ (\ie, $v$ is a node in $T$), then
$T^v$ denotes the fringe subtree of $T$ at $v$, \ie, the subtree consisting
of $v$ and all its descendants; $T^v$ is defined as a rooted tree with root
$v$.

$\Exp(\gl)$ denotes an exponential random variable with rate $\gl$; it thus
has mean $1/\gl$ and density function $\gl e^{-\gl x}$, $x>0$.

\subsection{Family trees of branching processes}\label{SSCMJ}

We follow \cite[Section 5]{SJ306}, to which we refer for further details. 
 Let $\cT_t$ be the family tree
of all individuals born up to time $t\ge0$ in a given 
Crump--Mode--Jagers (CMJ) process, starting at time $t=0$
with a single individual (the root).
Let the children of the root be born at (random) times $(\xi_i)_1^N$, where
$0\le N\le\infty$ and $0<\xi_1\le\xi_2\le\dots$.
We regard the (multi)set of birth times as a point process $\Xi$;
formally, $\Xi$ is the random (discrete) measure $\sum_i\gd_{\xi_i}$,
where $\gd_t$ is the  Dirac measure (point mass) at $t$.
Moreover, each individual $x$ has its own copy $\Xi_x$ of $\Xi$; 
the processes $\Xi_x$ are \iid{} (independent and identically distributed).
Let $\gs_x$ be the time individual $x$ is born.
For simplicity we assume that all individuals live forever.

Let $Z_t$ be the number of individuals at time $t$.
In the simplest, and most common, case, we define the stopping time
\begin{align}  \label{tau1}
\tau(n):=\inf\set{t:Z_t\ge n},
\end{align} 
the first time the number of individuals is at
least $n$, and $T_n:=\cT_{\tau(n)}$, the family tree at that time.
(By the assumptions below, $\tau(n)<\infty$ a.s.)
Thus $T_n$ is a random tree with $|T_n|\ge n$. Typically, the
birth times $\xi_i$ are continuous random variables and \as{} no two births
are simultaneous, and then $|T_n|=n$.

More generally, we fix a \emph{weight} $\psi(s)$.
This is assumed to be 
a \emph{characteristic}, i.e.,
a random function $\psi(s)\ge0$
associated to the root and its point process $\Xi$, and we assume that 
each individual $x$ is equipped with its own copy $\psi_x(s)$ of $\psi$;
 the simplest case is that $\psi_x(s)$ is a deterministic function of the
point process $\Xi_x$. (More generally, $\psi_x$ may also depend on the
entire tree of descendants of $x$, and possibly also on some extra randomness,
see \cite[in particular Remark 5.10]{SJ306}.)
We assume $\psi_x\in D\ooo$, and we exclude the trivial case $\psi(t)=0$ for
all $t\ge0$ a.s.
The argument $s\ge0$ of $\psi_x(s)$ should be
interpreted as the current age of $x$, which is $t-\gs_x$ at time $t$.
Let
\begin{align}
  Z_t^\psi:=\sum_{x:\gs_x\le t}\psi_x(t-\gs_x)
\end{align}
be the total weight at time $t\ge0$.
We then let
\begin{align}  \label{taupsi}
\tau(n):=\inf\set{t:Z_t^\psi\ge n},
\end{align} 
the first time the total weight is at
least $n$. (We define $\inf\emptyset=\infty$.)
Finally, as before, we define $T_n:=\cT_{\tau(n)}$.
Note that the choice $\psi(s)=1$ gives $Z_t^\psi=Z_t$, and thus the simple
definition \eqref{tau1} of $T_n$.

Examples of common random trees that can be constructed as $T_n$ in this way
are given in \refSs{SRRT}--\ref{Slast};
see further \cite{SJ306}.

Let $\mu:=\E\Xi$ be the intensity of the point process $\Xi$. 
In other words, $\mu$ is the
(deterministic) measure on $\ooo$ such that, for any Borel set $A$, $\mu(A)$ is
the expected number of children of the root born at times $t\in A$.
In particular, with $N\le\infty$ as above the (random) total number of
children of the root, $\mu\ooo=\E N\le\infty$.

We use the following assumptions throughout the paper:
\begin{PXenumerate}{A}
\item \label{BPfirst}
$\xi_1>0$, i.e., no children are born immediately at their parent's birth.
(Equivalently, $\mu\set0=0$.)
\item \label{BPnonlattice}
$\mu$ is not concentrated on any lattice $h\bbZ$, $h>0$.
(The results extend to the lattice case with suitable modifications, but we
do not know any interesting examples and ignore this case.) 
\item \label{BPsuper}
$N\ge1$ \as{} and $\E N>1$. 
(Thus, every individual has at least one child, so the process never
dies out, and $Z_\infty=\infty$ a.s.)
\item \label{BPmalthus}
There exists a real number $\ga$ (the \emph{Malthusian parameter}) such that
\begin{equation}\label{malthus}
\intoo e^{-\ga t}\mu(\ddx t) =1.  
\end{equation}
(By \ref{BPsuper}, $\ga>0$.)

\item \label{BPmub}
There exists $\gth<\ga$ such that
\begin{equation}\label{malthus-}
\intoo e^{-\gth t}\mu(\ddx t) <\infty.  
\end{equation}

%
\end{PXenumerate}

\subsection{Independence numbers}\label{SSind}
We collect here some simple and well-known properties of independence numbers
of (rooted) trees; see \eg{} \cite{MM1973}.

For a tree $T$, let $I(T)$ be the independence number of $T$, 
\ie, the maximum size of an independent set of nodes.
For a rooted tree $T$, let further $I_1(T)$ be the 
maximum size of an independent node set containing the root, 
and let $I_0(T)$ be the 
maximum size of an independent node set not containing the root.
Thus
\begin{align}\label{I01}
  I(T)=I_1(T)\lor I_0(T).
\end{align}
Furthermore, if the children of the root are $v_1,\dots,v_d$, then
it is easily seen that
\begin{align}
  I_0(T)&=\sum_{i=1}^d I(T^{v_i}),\label{I0}
\\
  I_1(T)&=1+\sum_{i=1}^d I_0(T^{v_i}).\label{I1}
\end{align}
Since $I_0\le I$ by \eqref{I01}, it follows that $I_1(T)\le I_0(T)+1$, and
thus
\begin{align}\label{I0I}
  I_0(T)\le I(T)\le I_0(T)+1.
\end{align}
Define
\begin{align}\label{iota}
  \iota(T):=I(T)-I_0(T)\in\setoi.
\end{align}
Then \eqref{I0} yields
\begin{align}
  I(T)=\iota(T)+\sum_{i=1}^d I(T^{v_i}),\label{Itoll}
\end{align}
which shows that the independence number $I(T)$ is an additive functional
on  rooted trees with toll function $\iota(T)$.

As is well known, \eqref{Itoll} is equivalent to
\begin{align}\label{Itoll1}
  I(T)=\sum_{v\in T} \iota(T^v).
\end{align}

Furthermore, by \eqref{iota}, \eqref{I01} and \eqref{I0}--\eqref{I1},
\begin{align}\label{iota+}
  \iota(T)=1
&\iff
I_1(T)=1+I_0(T)
\iff
\sumid I_0(T^{v_i}) = \sumid I(T^{v_i})
\notag\\&
\iff \iota(T^{v_i})=0 \text{ for every child $v_i$ of the root}.
\end{align}

Say that a node $v\in T$ is \emph{essential} if it belongs to every maximum
independent set of $T^v$. 
This is equivalent to $I_1(T^v)>I_0(T^v)$, and thus to $\iota(T^v)=1$.
In other words,
\begin{align}\label{ess}
  \iota(T^v)=\indic{v \text{ is essential}}.
\end{align}
In particular, $\iota(T)$ equals the indicator that the root is essential in
$T$. Note also that,
by \eqref{ess} and \eqref{iota+}, 
\begin{align}\label{noess}
  \text{a node is essential if and only if
none of its children is.}
\end{align}
\begin{remark}
By \eqref{Itoll1} and \eqref{ess},  the independence number $I(T)$ equals
the number of essential nodes in $T$.
Moreover, \eqref{noess} implies  that the set of essential nodes  is
independent, and thus an independent set of maximum size. 
\end{remark}

\section{Main result}

We next state our main theorem.
Recall that the Laplace functional of the point process $\Xi$ 
is defined as
\begin{align}\label{laplace}
  \LXi(f):=\E e^{-\int f\dd\Xi}
\end{align}
for (measurable) functions $f\ge0$ on $\ooo$.

\begin{theorem}\label{T1}
  Let $T_n$, $n\ge1$, be random trees that can be defined as stopped
family trees of   Crump--Mode--Jagers processes as in \refS{SSCMJ}, for some 
point process\/ $\Xi=(\xi_i)$ 
satisfying assumptions \AAA
and some weight $\psi$. Then, as \ntoo,
\begin{align}\label{t1}
  \frac{I(T_n)}{|T_n|}\asto \nu :=\ga\intoo e^{-\ga t} p(t)\dd t,
\end{align}
where $\ga$ is the Malthusian parameter and 
$p(t)$ is the unique function $\ooo\to\ooi$ satisfying
\begin{align}
  \label{tl1}
  p(t) 
&= \E \prod_{i:\xi_i\le t}\bigpar{1-p(t-\xi_i)}
= \E e^{\int_0^t\log(1-p(t-s))\dd\Xi(s)}
\notag\\&
= \LXi\Bigpar{-\log\bigpar{1-p(t-\cdot)}\indic{\cdot\le t}},
\qquad t\ge0
.
\end{align}
\end{theorem}

Note that the result does not depend on the choice of weigth $\psi$.

\begin{proof}
By \eqref{Itoll1} and \eqref{ess}, $I(T_n)/|T_n|$ is the fraction of nodes
in $T_n$ that are essential.
We apply \cite[Theorem 5.14(ii)]{SJ306}
to the property that a node is essential.
(This theorem is a special case of deep
results by Jagers and Nerman \cite{JagersNerman1984,NermanJagers1984},
see also \citet{Aldous}.)
Then \cite[(5.23)--(5.24)]{SJ306}
yield \eqref{t1}, with
\begin{align}\label{pt}
  p(t): 
=\P\bigpar{\text{the root is essential in $\ctt$}}
=\P\bigpar{\iota(\ctt)=1}
.\end{align}
To see \eqref{tl1}, 
condition on $\Xi=(\xi_i)_i$, \ie, on the
sequence of times that the root gives birth.
Then, the children of the root of $\ctt$ 
are the individuals $i$ born at times $\xi_i\le t$. Each such child has
grown a tree $\ctt^i$ that has the same distribution as $\cT_{t-\xi_i}$,
and thus
\begin{align}\label{alfa}
  \P\bigpar{{\iota(\ctt^i)=0}\mid\Xi}
=1-p(t-\xi_i),
\qquad \xi_i\le t.
\end{align}
Furthermore, still conditioned on $\Xi$, the events in \eqref{alfa} for
different $i$ are independent, and thus
using \eqref{iota+}, 
\begin{align}
  \P\bigpar{{\iota(\ctt)=1}\mid\Xi}
&=
  \P\bigpar{\iota(\ctt^i)=0 \text{ for every child $i$ of the root}\mid\Xi}
\notag\\&
=\prod_{i:\xi_i\le t}\bigpar{1-p(t-\xi_i)}.
\end{align}
Hence, \eqref{tl1} follows, using \eqref{laplace}.

Finally, suppose that $p_1(t)$ is another function $\ooo\to\oi$ that satisfies
\eqref{tl1}, and let $\gdp(t):=|p(t)-p_1(t)|$.
Then,
\begin{align}\label{gdp}
  \gdp(t)&
=\Bigabs{\E\prod_{i:\xi_i\le t}\bigpar{1-p(t-\xi_i)}
  -\E\prod_{i:\xi_i\le t}\bigpar{1-p_1(t-\xi_i)}}
\notag\\&
\le \E\sum_{i:\xi_i\le t}\bigabs{p(t-\xi_i)-p_1(t-\xi_i)}
= \E\sum_{i:\xi_i\le t}\gdp(t-\xi_i).
\end{align}
Fix $\gb>\ga$, and define $h(t):=\sup_{s\le t} e^{-\gb s}\gdp(s)\in\oi$.
Then, \eqref{gdp} yields
\begin{align}\label{qb}
    e^{-\gb t}\gdp(t) 
&\le{ e^{-\gb t}\E \sum_{i:\xi_i\le t} e^{\gb(t-\xi_i)}h(t)}
=h(t) \E  \sum_{i:\xi_i\le t} e^{-\gb\xi_i}
\notag\\&
= h(t)\E\intot e^{-\gb x}\dd\Xi(x)
= h(t)\intot e^{-\gb x}\dd\mu(x)
.\end{align}
Since $h(t)$ is monotone, this implies
\begin{align}\label{qj}
h(t)=\sup_{s\le t} \bigpar{   e^{-\gb s}\gdp(s) }
\le h(t)\intot e^{-\gb x}\dd\mu(x).
\end{align}
However, by \eqref{malthus},
\begin{align}
  \intot e^{-\gb x}\dd\mu(x)
\le   \intoo e^{-\gb x}\dd\mu(x)
<
\intoo e^{-\ga x}\dd\mu(x)=1,
\end{align}
and thus \eqref{qj} implies $h(t)=0$ for any $t\ge0$. Thus $p_1(t)=p(t)$,
and the solution to 
\eqref{tl1} is unique.
\end{proof}

Note that $\cT_0$ consists of the root only, and thus \eqref{pt} yields
the initial condition, also a trivial special case of \eqref{tl1},
\begin{align}\label{p0}
  p(0)=1.
\end{align}

\begin{remark}
  \label{Rcttau}
An explanation for the formula \eqref{t1} for $\nu$ is that a random fringe
tree
of $T_n$ converges in distribution to $\cttau$, 
the tree obtained by stopping the branching process 
at a time $\tau\sim\Exp(\ga)$ independent of the brancing process;
thus $\nu$ is the probability that the root of $\cttau$ is essential.
See further \cite{SJ306}.
\end{remark}

\begin{remark}\label{Rp<0}
  It is sometimes convenient to define $p(t):=0$ for $t<0$; then \eqref{tl1}
may be written 
\begin{align}\label{rp<0}
  p(t)=\E \prod_{i=1}^N\bigpar{1-p(t-\xi_i)},
\qquad t\ge0.
\end{align}
taking the product over all children (born yet or not). We will use this a
couple of times, but note that all formulas for $p(t)$ assume $t\ge0$.
\end{remark}

\section{The random recursive tree}\label{SRRT}

The random recursive tree is an example of a random tree that can be
constructed as in \refS{SSCMJ}, taking $\Xi$ to be a Poisson process with
constant intensity 1 on $\ooo$ and the trivial weight $\psi(t)=1$, see
\cite[Example 6.1]{SJ306}. 
Thus \refT{T1} applies and shows
$I(T_n)/|T_n|\asto\nu$ as \ntoo.

To find the limit $\nu$, note first that
by the standard formula
\cite[Theorem 3.9]{LastPenrose}
for the Laplace functional of a rate 1 Poisson process
\begin{align}
  \LXi(f)=e^{-\int(1-e^{-f(s)})\dd s},
\end{align}
 \eqref{tl1} yields
\begin{align}\label{rrt}
  p(t)
=\exp\Bigpar{-\int_0^t p(t-s)\dd s}
=\exp\Bigpar{-\int_0^t p(u)\dd u}.
\end{align}
This can also be seen directly as follows.
The number of children of the root at time $t$ has the Poisson distribution
$\Po(t)$. Furthermore, a child born at time $s\le t$ 
has probability $p(t-s)$ of being essential at time $t$, and thus,
by the  independence properties of the branching process, 
the children of the root that are essential at time $t$ are born according
to a random thinning of the rate 1 Poisson process;
 this thinning is a Poisson process on $[0,t]$ with intensity $p(t-\cdot)$.
In particular, the number of children of the root of $\ctt$ that are
essential is $\Po(\gl(t))$ with $\gl(t)=\int_0^t p(t-s)\dd s$.
By  \eqref{noess}, the root is essential if and only if this
number is 0, 
which has probability $e^{-\gl(t)}$,
and \eqref{rrt} follows.

Since $p(t)\in\oi$, \eqref{rrt} implies that $p(t)$ is continuous, and by
induction infinitely differentiable.
Differentiating \eqref{rrt} yields
\begin{align}
  p'(t)=-\exp\Bigpar{-\int_0^t p(u)\dd u} p(t) = -p(t)^2, \qquad t>0,
\end{align}
and thus
\begin{align}
  \frac{\ddx}{\dd t}\frac{1}{p(t)}=-\frac{p'(t)}{p(t)^2}=1.
\end{align}
Consequently, by the initial condition $p(0)=1$ \eqref{p0},
\begin{align}\label{rrp}
  p(t)=\frac{1}{1+t}.
\end{align}

 The Malthusian parameter $\ga=1$, and thus by \eqref{t1} and \eqref{rrp}
\begin{align}\label{nurrt}
  \nu = \intoo \frac{e^{-t}}{t+1}\dd t
= \intoi \frac{1}{1-\log x}\dd x
= 0.59634736\dots
\end{align}
as proved by
\citet{MM1975,MM1978},
\citet{Dadedzi} and \citet{FuchsEtal}; 
this number is
known as the Euler--Gompertz constant.

\section{The Binary Search Tree}\label{SBST}

The binary search tree is another example where \refT{T1} applies. 
Now each node gets two children,
after waiting times that are independent and $\Exp(1)$.
Again, $\psi(t)=1$.

We proceed to find the limit $\nu$.
In this case, \eqref{tl1} yields
\begin{align}\label{pbst}
  p(t)=\Bigpar{1-\int_0^t p(t-s)e^{-s}\dd s}^2
=\Bigpar{1-e^{-t}\int_0^t e^up(u)\dd u}^2,
\qquad t\ge0.
\end{align}
This can also, perhaps more easily, be seen as follows.
As always, if a child is born at time $s\le t$, then the probability that
this child is essential at time $t$ is $p(t-s)$. Hence, the probability that
the left child of the root is born and is essential at time $t$ is
$\int_0^t p(t-s)e^{-s}\dd s$, and thus the probability that there is 
no left child that is essential equals
$1-\int_0^t p(t-s)e^{-s}\dd s$.
The same holds for the right child, and since the two children appear and
develop independently, \eqref{noess} yields \eqref{pbst}.

To solve the functional equation \eqref{pbst},
let $g(t):=p(t)\qq$, so \eqref{pbst} may be written
\begin{align}\label{gbst}
  g(t)
={1-e^{-t}\int_0^t e^sg(s)^2\dd s}.
\end{align}
It follows, by induction, that $g$ is infinitely differentiable; 
furthermore, \eqref{gbst} yields the differential equation
\begin{align}\label{gbstt}
  g'(t)
=e^{-t}\int_0^t e^sg(s)^2\dd s -g(t)^2
=1-g(t)-g(t)^2.
\end{align}
This differential equation is separable and can be written
\begin{align}\label{dgt}
  \frac{\ddx g}{1-g-g^2}={\ddx t},
\end{align}
which is solved by standard methods as follows.

Let $\gam_\pm:=\frac{-1\pm\sqrt5}{2}$ be the roots of $1-\gam-\gam^2=0$.
Then
\begin{align}
  \frac{1}{1-g-g^2}
=-\frac{1}{(g-\gam_+)(g-\gam_-)}
=\frac{1}{\gam_+-\gam_-}\Bigpar{\frac{1}{g-\gam_-}-\frac{1}{g-\gam_+}}
\end{align}
and thus \eqref{dgt} can be integrated to
\begin{align}\label{aba}
 \log(g(t)-\gam_-)-\log(g(t)-\gam_+)=(\gam_+-\gam_-)t+C,
\end{align}
where $g(0)=1$ yields
$C= \log(1-\gam_-)-\log(1-\gam_+)$.
Hence,
\begin{align}\label{abb}
  \frac{g(t)-\gam_-}{g(t)-\gam_+}
=\frac{1-\gam_-}{1-\gam_+} e^{(\gam_+-\gam_-)t}
\end{align}
and thus
\begin{align}\label{olle}
  g(t)&=
\frac{\gam_+(1-\gam_-)e^{(\gam_+-\gam_-)t}-\gam_-(1-\gam_+)}
{(1-\gam_-)e^{(\gam_+-\gam_-)t}-(1-\gam_+)}
\notag\\&=
\frac{\gam_+(1-\gam_-)e^{\gam_+t}-\gam_-(1-\gam_+)e^{\gam_-t}}
{(1-\gam_-)e^{\gam_+t}-(1-\gam_+)e^{\gam_-t}}.
\end{align}
Consequently, $p(t)=g(t)^2$ with $g(t)$ given by \eqref{olle}.
Note also that we have
$\gam_+-\gam_-=\sqrt5$ and,
with $\phi:=\frac{1+\sqrt5}2$, the golden ratio,
\begin{align}\label{gap}
  \gam_+ &=\phi\qw,&
1-\gam_+&=\phi\qww,
\\\label{gam}
\gam_-&=-\phi,&
1-\gam_-&=\phi^2.
\end{align}

The Malthusian parameter $\ga=1$, and thus \eqref{t1} and \eqref{olle}
yield,
with $x=e^{-t}$,
\begin{align}\label{pelle}
  \nu
&=
\intoo e^{-t}g(t)^2\dd t 
=\intoo e^{-t}\biggparfrac{\phi e^{\sqrt5\,t}+\phi\qw}
{\phi^2e^{\sqrt5\, t}-\phi\qww}^2\dd t
\notag\\&=
\intoi\biggparfrac{\phi +\phi\qw x^{\sqrt5}}
{\phi^2-\phi\qww x^{\sqrt5}}^2\dd x
=
\frac{\phi^2}{\sqrt5}
\intoi\biggparfrac{\phi^2+y}
{\phi^4- y}^2y^{1/\sqrt5-1}\dd y.
\end{align}
This integral can be evaluated as the sum of
a rapidly (geometrically) convergent
series by expanding 
$(\phi^4-y)^{-2}=\phi^{-8}(1-\phi^{-4}y)^{-2}$
into a power series,  
which yields
\begin{align}\label{kalle}
  \nu
=
\sumko(k+1)\phi^{-4k-6}
\Bigpar{\frac{\phi^4}{k\sqrt5+1} + \frac{2\phi^2}{(k+1)\sqrt5+1}
+\frac{1}{(k+2)\sqrt5+1}}.
\end{align}
The integrals \eqref{pelle} and the sum \eqref{kalle} are all easily
evaluated numerically, yielding (by Maple)
\begin{align}
  \nu= 0.54287631\dots   
\end{align}
as found by 
\cite{Dadedzi} and \cite{FuchsEtal}.

\section{Preferential attachment trees}\label{SPA}

Consider now a preferential attachment tree, where nodes are added one by
one, and each new node chooses a parent at random, with the probability of
choosing a node $v$ as the parent is proportional to $\chi d(v)+\rho$, where
$d(v)$ is the current outdegree of $v$, and $\chi$ and $\rho$ are 
given constants.
(Here $\rho>0$ and either $\chi\ge0$ or $\rho/|\chi|$ is an integer. Only
the ratio $\chi/\rho$ is significant.)
This random tree can be constructed by a CMJ process where an individual
that already has $k$ children gets the next child with rate $\chi k+\rho$;
see
\cite[Example 6.4]{SJ306}. Again the weight $\psi(t)=1$.
Thus \refT{T1} applies and shows
$I(T_n)/|T_n|\asto\nu$ as \ntoo.

To find $p(t)$ and $\nu$, 
we use instead of \eqref{tl1} the following (closely related) argument.
Consider also,  for  $\gl>0$, 
a modified branching process $\fX_\gl$,
where the
starting individual (= the root) is special, and gets children with the rate 
$\chi k+\gl$, where $k$ is the current number  of children. All other
individuals are as before, with rate $\chi k+\rho$.
(If $\chi<0$, we assume that $\gl/|\chi|$ is an integer; this case will be
enough below.)
Let $p_\gl(t)$ be the probability that the root is essential in the 
family tree of $\fX_\gl$  at time $t$. 
Note that if $\gl=\rho$, then the modified process equals the original one,
and thus $p(t)=p_\rho(t)$.

Consider again the original process, let $t>0$, 
and condition on the first
child of the root being born at time $s\le t$. 
The probability that this child is not essential at time $t$ is $1-p(t-s)$.
Furthermore, if we ignore this child and its descendants, the rest of the
tree
evolves after time $s$ as the modified process $\fX_{\chi+\rho}$. 
Hence, the probability that  no child of the root after
the first is essential at time $t$ is $p_{\chi+\rho}(t-s)$.
Consequently, conditioned on the first child being born at time $s\le t$,
the probability that the root has no essential child at time $t$ is
$\bigpar{1-p(t-s)}p_{\chi+\rho}(t-s)$. By \eqref{noess}, this is the
conditional probability that the root is essential at time $t$, and since
the time the first child is born has the distribution $\Exp(\rho)$, we have
\begin{align}\label{bla}
  p(t)&
=e^{-\rho t} +\rho\intot e^{-\rho s}\bigpar{1-p(t-s)}p_{\chi+\rho}(t-s)\dd s
\notag\\&
=e^{-\rho t} +\rho e^{-\rho t}\intot e^{\rho u}\bigpar{1-p(u)}p_{\chi+\rho}(u)\dd u.
\end{align}

If we have two independent modified processes  $\fX_{\gl_1}$
and $\fX'_{\gl_2}$, then we may merge them by identifying the two roots. 
This
yields a modified process $\fX_{\gl_1+\gl_2}$ with parameter $\gl=\gl_1+\gl_2$.
 There are no essential children of the root in the combined process if
and only if there are none in both modified processes taken separately;
hence,
it follows that, for any $t\ge0$ and $\gl_1,\gl_2>0$,
\begin{align}\label{putte}
  p_{\gl_1+\gl_2}(t)=  p_{\gl_1}(t)  p_{\gl_2}(t).
\end{align}
Fix $t>0$. Since $0< p_{\gl}(t)\le 1$, \eqref{putte} implies that $\gl\mapsto
p_\gl(t)$ is decreasing, which in turn implies that
\eqref{putte} has the solution $p_{\gl}(t)=e^{-C(t)\gl}$ for some $C(t)\ge0$.
Hence, 
\begin{align}\label{gron}
  p_\gl(t)=p_\rho(t)^{\gl/\rho}
=p(t)^{\gl/\rho}.
\end{align}
Combining \eqref{bla} and \eqref{gron} yields the functional equation
\begin{align}\label{brun}
  p(t)
=e^{-\rho t} +\rho e^{-\rho t}\intot e^{\rho u}\bigpar{1-p(u)}p(u)^{\chi/\rho+1}\dd u.
\end{align}
Again, $p$ is infinitely differentiable, and
taking the derivative yields
\begin{align}\label{p'}
  p'(t)&
=-\rho e^{-\rho t} 
-\rho^2 e^{-\rho t}\intot e^{\rho u}\bigpar{1-p(u)}p(u)^{\chi/\rho+1}\dd u
+\rho \bigpar{1-p(t)}p(t)^{\chi/\rho+1}
\notag\\&
=
-\rho p(t)
+\rho \bigpar{1-p(t)}p(t)^{\chi/\rho+1}.
\end{align}
Let $\chix:=\chi/\rho$, and 
\begin{align}\label{h}
  h(x):=
x-x^{\chix+1}(1-x)
=x-x^{\chix+1}+x^{\chix+2}.
\end{align}
Then \eqref{p'} can be written
\begin{align}\label{pip}
  p'(t)=-\rho h\xpar{p(t)}.
\end{align}

By \eqref{h}, $h(1)=1>0$ and $h(0)=0$.
(If $\chi<0$, then $\rho=m|\chi|$ for an integer $m\ge2$, 
and thus $\chi'=-1/m\in[-\frac12,0)$, so $\chi'+1>0$ also in this case.)
Let $q$ be the largest zero of $h$ in $\oi$, \ie,
\begin{align}
  q:=\max\set{x\in\oi:h(x)=0}.
\end{align}
By continuity and $h(0)=0$, this maximum always exists. 
Furthermore, if $\chi\ge0$, then \eqref{h} implies $h(x)>0$ on $(0,1]$, 
and thus $q=0$.
On the other hand, if $\chi<0$, 
then $h(x)<0$ for small positive $x$,  
and thus $0<q<1$.

The function $p(t)$ is continuous on $\ooo$, with $p(0)=1$.
Suppose that $p(t)=q$ for some $t<\infty$, and let $t_0$ be the smallest 
such $t$. Then $p(t)\in[q,1]$ for $t\in[0,t_0]$. However, $h(x)$ is
continuously differentiable and thus Lipschitz on $[q,1]$, as is seen by
considering the cases $\chi\ge0$ and $\chi<0$ separately, and 
thus the differential equation \eqref{pip} has at most one solution for
$t\in[0,t_0]$ with $p(t)\in[q,1]$ and $p(t_0)=q$.
Since $p(t)=q$
is another solution of \eqref{pip}, this is a contradiction.
Hence, $p(t)\neq q$, and thus, by continuity, $p(t)>q$ for all $t\ge0$.

This further implies, by \eqref{pip} again,  
that $p(t)$ is strictly decreasing on $\ooo$.
Hence the limit $p(\infty):=\lim_{\ttoo}p(t)$ exists.
Then \eqref{pip} implies $p'(t)\to -\rho h(p(\infty))$ as \ttoo,
and thus $p(\infty)>q$ is impossible; hence,
\begin{align}
  p(t)\to p(\infty)=q,
\qquad\ttoo.
\end{align}
Consequently, $p$ is a bijection $\ooo\to(q,1]$.

Let, for $x\in (q,1]$,
\begin{align}\label{Psi}
  \Psi(x):=\int_x^1\frac{1}{h(y)}\dd y.
\end{align}
Thus $\Psi(1)=0$ and $\Psi'(x)=-1/h(x)$.
Hence, \eqref{pip} implies by the chain rule
\begin{align}
\frac{\ddx}{\dd t}  \Psi(p(t))=
\Psi'(p(t))p'(t)= \frac{-\rho h(p(t))}{-h(p(t))}
=\rho,
\end{align}
and thus 
\begin{align}\label{abba}
  \Psi(p(t))=\rho t,
\qquad t\ge0.
\end{align}
Hence, if $\Psi\qw:\ooo\to(q,1]$ denotes the inverse function, then
\begin{align}\label{pipp}
  p(t)=\Psi\qw(\rho t).
\end{align}
It follows from \eqref{abba}, letting \ttoo, that
the limit $\Psi(q)=\infty$, which also easily can be seen directly from
\eqref{Psi}.

The Malthusian parameter $\ga=\chi+\rho$, see \cite[(6.20)]{SJ306},
and thus \eqref{t1} and \eqref{pipp} yield
\begin{align}
  \nu = \ga\intoo e^{-\ga t} \Psi\qw(\rho t)\dd t
=(\chix+1)\intoo e^{-(\chix+1)s}\Psi\qw(s)\dd s.
\end{align}
The change of variables $s=\Psi(x)$ 
and an integration by parts
yield the formulas
\begin{align}\label{nugpa1}
  \nu &=
(\chix+1) \int_{1}^q e^{-(\chix+1)\Psi(x)} x \Psi'(x)\dd x
\\& \label{nugpa2}
=
(\chix+1) \int_{q}^1 e^{-(\chix+1)\Psi(x)} \frac{x}{h(x)} \dd x
\\& \label{nugpa3}
=1- \int_{q}^1 e^{-(\chix+1)\Psi(x)} \dd x.
\end{align}
These integrals can be evaluated numerically.

\begin{example}
  Let $\chi=0$ and $\rho=1$; this yields the random recursive tree in
  \refS{SRRT}.
We have $\chix=0$, and thus \eqref{h} yields $h(x)=x^2$.
Hence, \eqref{Psi} yields $\Psi(x)=1/x-1=(1-x)/x$,
and \eqref{pipp} yields  \eqref{rrp} again.
Furthermore, \eqref{nugpa2} becomes
\begin{align}\label{nurrt-g}
  \nu = \intoi e^{1-1/x} x\qw \dd x
= \int_1^\infty e^{1-y} y\qw \dd y,
\end{align}
which by a change of variables agrees with \eqref{nurrt}.
\end{example}

\begin{example}
  Let $\chi=\rho=1$; this yields the standard preferential attachment random
  tree. 
(This is the same as the plane oriented recursive tree 
 \cite{Szymanski}; it is a special case of the preferential attachment
 graphs \cite{BarabasiA}, \cite[Chapter~8]{HofstadI}.)
We have $\chix=1$, and thus \eqref{h} yields 
\begin{align}
  h(x)=x-x^2+x^3.
\end{align}
We find from \eqref{Psi},
\begin{align}
  \Psi(x)=\tfrac12\log(x^2-x+1)-\log(x)-\frac{1}{\sqrt3} 
\arctan \Bigparfrac{2x-1}{\sqrt3}
+\frac{\pi}{6\sqrt 3}
\end{align}
and thus \eqref{nugpa3} yields, with an integral that magically 
has an elementary primitive function,
\begin{align}
  \nu &= 
1-\intoi e^{-2\Psi(x)}\dd x
= 
1-e^{-\pi/3\sqrt3}\intoi e^{\frac{2}{\sqrt3}\arctan((2x-1)/\sqrt3)}\frac{x^2}{1-x+x^2}\dd x
\notag\\&
= 
1-e^{-\pi/3\sqrt3}\Bigsqpar{(x-1) e^{\frac{2}{\sqrt3}\arctan((2x-1)/\sqrt3)}}_0^1
\notag\\&= 
1-e^{-2\pi/3\sqrt3}
= 0.70156394\dots
\end{align}
Thus, \refT{T1}  shows that for the standard preferential attachment tree,
$I(T_n)/|T_n|\asto 1-e^{-2\pi/3\sqrt3}$ as \ntoo.
\end{example}

\section{Extended binary search trees}\label{SXBST}

An extended binary search tree is a binary search tree where we have added
further leaves at all possible places;
thus the original nodes (called internal nodes) have all two children each,
and the new nodes (called external nodes) have no children.
This can be constructed by a CMJ process where each individual gets twins
after an $\Exp(1)$ time (and no further children). 
Note that in the tree $\ctt$, the internal nodes
are the ones that have had children, while the others are external nodes.

We may choose to measure the size of an extended binary search tree in
three different ways: the total number of nodes, the number of internal
nodes, or the number of external nodes. (These are related in simple ways,
since in the binary case treated here, the number of external nodes is
always 1 + the number of internal nodes.)
We obtain these three versions as our $T_n$ 
by choosing different weight functions
$\psi$;
$\psi(t)=1$ as usual gives the total number of vertices, while
the number of internal vertices is given by $Z_t^\psi$ with
$\psi(t):=\indic{\xi_1\le t}$ and the number of external vertices 
is given by $\psi(t):=\indic{\xi_1> t}$.
Recall that \refT{T1} applies, and gives the same limit, to all three
versions.

Since we have $\xi_1=\xi_2$, \eqref{tl1} yields
(with $p(t)=0$ for $t<0$, see \refR{Rp<0})
\begin{align}
  p(t)&
=\E\bigsqpar{ \bigpar{1-p(t-\xi_1)}^2}
=e^{-t}+\intot e^{-s} \bigpar{1-p(t-s)}^2\dd s
\notag\\&
=e^{-t}+e^{-t}\intot e^{u} \bigpar{1-p(u)}^2\dd u,
\qquad t\ge0.
\end{align}
This yields the differential equation
\begin{align}\label{pbex}
  p'(t) = -p(t) + \bigpar{1-p(t)}^2
= 1-3p(t)+p(t)^2.
\end{align}
Let $q(t):=1-p(t)$. Then \eqref{pbex} yields
\begin{align}
  q'(t)=-p'(t)=1-q(t)-q(t)^2,
\end{align}
which is the same differential equation as \eqref{gbstt}, although now the
initial condition is $q(0)=1-p(0)=0$.
The general solution is as in \eqref{aba}, and we obtain, \cf{}
\eqref{abb}--\eqref{olle},
\begin{align}\label{abc}
  \frac{q(t)-\gam_-}{q(t)-\gam_+}
=\frac{\gam_-}{\gam_+} e^{(\gam_+-\gam_-)t}
\end{align}
and, using \eqref{gap}--\eqref{gam},
\begin{align}
  q(t)=\frac{e^{(\gam_+-\gam_-)t}-1}{\gam_+-\gam_-e^{(\gam_+-\gam_-)t}}
=\frac{e^{\sqrt5 t}-1}{\phi\qw+\phi e^{\sqrt 5t}}
=\frac{\phi e^{\sqrt5 t}-\phi}{\phi^2 e^{\sqrt 5t}+1}
.\end{align}
Thus, recalling $\phi^2=\phi+1$,
\begin{align}
  p(t)=1-q(t)=
\frac{(\phi^2-\phi)e^{\sqrt5 t}+\phi+1}{\phi^2 e^{\sqrt 5t}+1}
=\frac{ e^{\sqrt5 t}+\phi^2}{\phi^2 e^{\sqrt 5t}+1}
\end{align}
and, since the Malthusian parameter $\ga=1$, using $x=e^{-t}$,
\begin{align}\label{nuxbst}
  \nu&
=\intoo p(t)e^{-t}\dd t
=\intoi \frac{1+\phi^2 x^{\sqrt5}}{\phi^2+x^{\sqrt 5}}\dd x
=\phi^2-\xpar{\phi^4-1}\intoi \frac{\ddx x}{\phi^2+ x^{\sqrt5}}
\notag\\&
=\phi^2-\xpar{3\phi+1}\sumko (-1)^k \frac{\phi^{-2-2k}}{k\sqrt 5+1}
=  0.5987899\dots
\end{align}

\section{$m$-ary search trees}\label{Smary}

Consider now $m$-ary search tree, for a given $m\ge3$. (The case $m=2$ was
studied in \refS{SBST}.)
The $m$-ary search tree $T_n$ generated by $n$ random keys
can be constructed by the following CMJ process and weight $\psi$, see
\cite[Section 7.2]{SJ306}.

Each individual (node) starts by gaining weight; the weight $\psi(t)$
represents the 
number of keys in the node. It starts with $\psi(0)=1$, and then increases
by 1 after successive independent waiting times $Y_2,\dots,Y_{m-1}$ with
$Y_i\sim\Exp(i)$. At time $S:=\sum_{i=2}^{m-1}Y_i$ the weight thus reaches $m-1$;
this marks puberty, and the node becomes fertile and gets $m$ children
after further independent waiting times $X_i\sim\Exp(1)$. (Thus, child $i$
is born at $S+X_i$.)

\refT{T1} thus applies. To find $\nu$, we 
condition on $S$ and find 
that if
$0\le s\le t$, then
(with $p(u)=0$ for $u<0$, see \refR{Rp<0})
\begin{align}
  \E&\Bigpar{\prod_{i:\xi_i\le t}\bigpar{1-p(t-\xi_i)}\mid S=s}
=
  \E\Bigpar{\prod_{i=1}^m\bigpar{1-p(t-\xi_i)}\mid S=s}
\notag\\&
= \E{\prodim\bigpar{1-p(t-s-X_i)}}
=\Bigpar{1-\int_0^{t-s}e^{-x}p(t-s-x)\dd x}^m
\notag\\&
=\Bigpar{1-e^{s-t}\int_0^{t-s}e^{y}p(y)\dd y}^m
.\end{align}
Hence, \eqref{tl1} or \eqref{rp<0} yields, 
for $t\ge0$, if $f_S$ is the density function of $S$,
\begin{align}\label{pm}
p(t)=\P(S>t)+\int_{0}^t \Bigpar{1-e^{s-t}\int_0^{t-s}e^{y}p(y)\dd y}^mf_S(s)\dd s.
\end{align}
(We assume $m\ge3$; for $m=2$ we have $S=0$ and \eqref{pm} is replaced by
\eqref{pbst}.)  
We define
\begin{align}\label{gm}
  g(t):=1-e^{-t}\intot e^y p(y)\dd y
\end{align}
and write \eqref{pm} as
\begin{align}
  \label{pm2}
p(t)&=P(S>t)+\intot g(t-s)^mf_S(s)\dd s
\notag\\&
=P(S>t)+\intot g(s)^mf_S(t-s)\dd s
.\end{align}

For simplicity, we consider in the sequel only the case $m=3$.
Then $S=Y_2\sim \Exp(2)$, and \eqref{pm2} becomes
\begin{align}
  \label{pm3}
p(t)&
=e^{-2t}+ 2e^{-2t}\intot e^{2s}g(s)^3\dd s
.\end{align}
It follows from \eqref{gm} and \eqref{pm3} by induction that $g(t)$ and
$p(t)$ are infinitely differentiable on $\ooo$, and differentiation yields
\begin{align} 
  g'(t)&=1-g(t)-p(t),\label{g3}
\\
p'(t)&=-2p(t)+2g(t)^3.\label{p3}
\end{align}
These equations can (as far as we know) only be solved numerically.
and then $\nu$ can be computed numerically by \eqref{tl1}, with
the Malthusian parameter $\ga=1$ \cite{SJ306}.
We obtain
\begin{align}
  \nu=\intoo e^{-t}p(t)\dd t
     = 0.58705155\dots
\end{align}
(Actually, we consider the system 
\set{\eqref{p3},\;
\nu'(t)=e^{-t}p(t)} with $p(0)=1$ and $\nu(0)=0$, and use Maple to find
$\nu=\nu(\infty)$.)

\begin{remark}
  Extended $m$-ary search trees may be considered similarly as the case
  $m=2$ in \refS{SXBST}, see \cite[Section 7.1]{SJ306}.
We leave this case to the reader.
\end{remark}

\label{Slast}

\begin{ack}
  I thank Cecilia Holmgren and Stephan Wagner for helpful discussions and
  comments. 
\end{ack}

\newcommand\AAP{\emph{Adv. Appl. Probab.} }
\newcommand\JAP{\emph{J. Appl. Probab.} }
\newcommand\JAMS{\emph{J. \AMS} }
\newcommand\MAMS{\emph{Memoirs \AMS} }
\newcommand\PAMS{\emph{Proc. \AMS} }
\newcommand\TAMS{\emph{Trans. \AMS} }
\newcommand\AnnMS{\emph{Ann. Math. Statist.} }
\newcommand\AnnPr{\emph{Ann. Probab.} }
\newcommand\CPC{\emph{Combin. Probab. Comput.} }
\newcommand\JMAA{\emph{J. Math. Anal. Appl.} }
\newcommand\RSA{\emph{Random Structures Algorithms} }
\newcommand\DMTCS{\jour{Discr. Math. Theor. Comput. Sci.} }

\newcommand\AMS{Amer. Math. Soc.}
\newcommand\Springer{Springer-Verlag}
\newcommand\Wiley{Wiley}

\newcommand\vol{\textbf}
\newcommand\jour{\emph}
\newcommand\book{\emph}
\newcommand\inbook{\emph}
\def\no#1#2,{\unskip#2, no. #1,} 
\newcommand\toappear{\unskip, to appear}

\newcommand\arxiv[1]{\texttt{arXiv:#1}}
\newcommand\arXiv{\arxiv}

\def\nobibitem#1\par{}

\end{document}